\documentclass[a4paper,english]{amsart}
\usepackage{latexsym}
\usepackage{amsmath}
\usepackage{babel}

\usepackage{amssymb}
\newcommand{\be}{\begin{equation}}
\newcommand{\ee}{\end{equation}}

\newtheorem{theorem}{Theorem}[section]

\newtheorem{lemma}[theorem]{Lemma}

\newtheorem{definition}[theorem]{Definition}
\newtheorem{definition and theorem}[theorem]{Definition and Theorem}
\newtheorem{remark}[theorem]{Remark}
\newtheorem{*remark}[theorem]{$^* $Remark}

\newtheorem{*exercise}[theorem]{$^* $Exercise}
\newtheorem{**exercise}[theorem]{$^{** } $Exercise}

\begin{document}
\selectlanguage{english}

\title[VaR and Expected Shortfall for Elliptic Linear Portfolios
]{ Value-at-Risk and Expected Shortfall For Linear Portfolios with
Elliptically Distributed Risk Factors}
\author[Jules Sadefo Kamdem]{Jules SADEFO KAMDEM \\
     \newline
Laboratoire de Math\'ematique \\
CNRS UMR 6056 \\
Universit\'e de Reims }
\thanks{This draft is a part  of J.SADEFO-KAMDEM (university of Reims)PhD Thesis.
It has been presented at the workshop on modelling and computation
in Financial Engineering at Bad Herrenalb, Germany May 6-8,2003.\\
Many thanks for the comments of Professor R.Brummelhuis.  \\
\underline{Author's address}: University of Reims, Laboratoire de
Math\'ematique UMR 6056-CNRS , BP 1039 Moulin de la
Housse , 51687 Reims cedex 2 FRANCE.\\
e-mail: sadefo@univ-reims.fr}
\maketitle
\begin{abstract}
   In this paper, we generalize the parametric $\Delta $-VaR
   method from portfolios with normally distributed risk factors to
   portfolios with elliptically distributed ones.
   We treat both the expected shortfall and the Value-at-Risk of
   such portfolios. Special attention is given to the particular case
   of a multivariate $t $-distribution.
\end{abstract}
\bigskip


\section{{\bf Introduction} }

 \medskip

The original RiskMetrics methodology for estimating VaR was based
on parametric methods, and used the multi-variate normal
distribution. This approach works well so-called {\em linear
portfolios}, that is, those portfolios whose aggregate return is,
to a good approximation, a linear function of the returns of the
individual assets which make up the portfolio, and in situations
where the latter can be assumed to be jointly normally
distributed. For other portfolios, like portfolios of derivatives
depending non-linearly on the return of the underlying, or
portfolios of non-normally distributed assets, one generally turns
to Monte Carlo methods to estimate the VaR. Monte Carlo
methodology has the obvious advantage of being almost universally
applicable, but has the disadvantage of being much slower than
comparable parametric methods, when the latter are available. This
is an issue in situations demanding for real-time evaluation of
financial risk. For non-linear portfolios, practitioners, as an
alternative to Monte Carlo, use $\Delta $-normal VaR methodology,
in which the portfolio return is linearly approximated, and an
assumption of normality is made. Such methods present us with a
trade-off between accuracy and speed, in the sense that they are
much faster than Monte-carlo, but much less accurate unless the
linear approximation is quite good, and the normality hypothesis
holds well. In case the linear approximation is of poor quality,
or inherently instable (as is for example the case when the
portfolio $\Delta $ is close to 0), one turns to a higher order
approximations, for example the quadratic one, while keeping the
normality assumption. This leads to the so-called $\Gamma - \Delta
$ VaR, which can be evaluated using Mont Carlo, but for which also
a number of semi-parametric methods have been developed (see for
example Duffy and Pan \cite{DP}, and others.

An obvious first generalization is to keep the linearity
assumption, but replace the normal distribution by some other
family of multi-variate distributions. This is the subject of the
present paper. As an alternative to the normal hypothesis we will
assume that the log returns of the portfolio's constituents have a
multivariate elliptic distribution in one of the elliptic classes
$N (\mu , \Sigma , \phi ) $, cf. section 1 below for the precise
definition. These have the advantage that,like normal
distributions, their dependence structure is completely determined
by their mean and variance, once a choice is made for $\phi $.
Particular examples are the normal distributions and the
multi-variant Student-$t $ distributions. The latter are an
obvious first choice since they posses heavy tails.
\medskip

Glasserman, Heidelberger and Shahabuddin \cite{GHS} present a method to
compute even $\Gamma - \Delta $ VaR using a semi-parametric method
based on the Fourier transform, but their methodology seems to be
restricted to $t $-distributions. See also Lopez and Walter \cite{LW}
and references therein for further applications of the $t $
distribution to VaR. Note that one shortcoming of the multivariate
t-distribution is that all the marginal distributions must have
the same degrees of freedom, which implies that every risk factor
has equally heavy tails.

In a sequel paper we intend to extend the present analysis to
include quadratically non-linear portfolios. We note that in a
companion paper we give very precise analytic VaR estimates for
the VaR of such quadratical portfolios, in case the underlying
risk-factors follow another class of non-Gaussian distributions,
namely the multivariate generalized normal distributions.

\medskip \quad The paper
is organized, as follows: In section 2, we will analyze the VaR of
a linear portfolio with elliptically distributed risk-factors,
paying special attention to the case of a multi-variate Student
distribution. Used in conjunction with a first order Taylor
approximation of a portfolio's Profit \& Loss function, this will
give rise to the notions of {\em Delta-Elliptic VaR} and {\em
Delta-Student VaR}, in analogy with the familiar Delta-Normal VaR.
We show, for example, how to reduce the computation of the
Delta-Student VaR to finding the zero's of a certain special
function. In section 3 we show how to extend our procedure to
mixtures of elliptic distributions. Section 4 treats the expected
shortfall for general elliptic linear portfolios and for the
special case of Student ones. Finally, in section 5 we discuss
some potential application areas.

\section{Linear Portfolio VaR with elliptic distributions}
\medskip

In this section we perform a parametric analysis that relies on
the assumption that the pricing function of the portfolio is
linear in the risk factors. Note that parametric methods provide
very fast answers which are, however, only as accurate as the
underlying linearity assumption.
\medskip

We will use the following notational conventions for the action of
matrices on vectors: single letters $x , y , \cdots $ will denote
{\em row vectors} $(x_1 , \cdots , x_n ) $, $(y_1 , \cdots y_n )
$. The corresponding column vectors will be denoted by $x^t , y^t
$,the $^t $ standing more generally for taking the transpose of
any matrix. Matrices $A = (A_{ij } )_{i, j } $, $B $ , etc. will
be multiplied in the usual way. In particular, $A $ will act on
vectors by left-multiplication on column vectors, $A y^t $, and by
right multiplication on row vectors, $xA $; $x \cdot x = x x^t =
x_1 ^2 + \cdots + x_n ^2 $ will stand for the Euclidean inner
product.
\medskip

A portfolio with time-$t $ value $\Pi (t) $ is called linear if
its profit and loss $\Delta \Pi(t)= \Pi (t) - \Pi (0) $ over a
time window, [0 t] is a linear function of the returns
$X_{1}(t),\ldots,X_{n}(t)$ of its constituents over the same time
period:
$$
\Delta \Pi(t)= \delta_{1} X_{1}+\delta_{2} X_{2}+...+\delta_{n}
X_{n}
$$
This will for instance be the case for ordinary portfolios of
common stock, if we use percentage returns, and will also hold to
good approximation with log-returns, provided the time window
[0,t] is small. We will drop the time $t $ from our notations,
since it will be kept fixed, and simply write $X_{j}$,$\Delta
\Pi$, etc. We also put
$$
X = (X_1 , \cdots , X_n ) ,
$$
so that $\Delta \Pi = \delta \cdot X = \delta X^t . $
\medskip

We now assume that the $X_j $ are elliptically distributed with
mean $\mu$ and correlation matrix $\Sigma=AA^t$:
$$
(X_{1},\ldots,X_{n})\sim N(\mu,\Sigma , \phi ).
$$
This means that the pdf of X is of the form
$$
f_{X}(x)=|\Sigma|^{-1} g((x-\mu)\Sigma^{-1}(x-\mu)^t ) ,
$$
where $|\Sigma | $ stands for the determinant of $\Sigma $, and
where $g: \mathbb{R }_{\geq 0 } \to 0 $ is such that the Fourier
transform of $g(|x|^2 ) $, as a generalized function on $\mathbb{R
}^n $, is equal to $\phi (|\xi |^2 ) $\footnote{One uses $\phi $
as a parameter for the class of elliptic distributions, since it
is always well-defined as a continuous function: $\phi (|\xi |^2 )
$ is simply the characteristic function of an $X \sim N(0, Id ,
\phi ) $. Note, however, that in applications we'd rather know $g
$}. Assuming that $g $ is continuous, and non-zero everywhere, the
Value at Risk at a confidence level of $1 - \alpha$ is given by
solution of the following equation: \be \mbox{Prob } \{ \Delta
\Pi(t) < - VaR_{\alpha}\}= \alpha    \label{1} \ee Here we follow
the usual convention of recording portfolio losses by negative
numbers, but stating the Value-at-Risk as a positive quantity of
money.

In terms of our elliptic distribution parameters we have to solve
the following equation:
$$
\alpha ={|\Sigma|}^{-1/2} \int_{\{\delta \cdot x  \leq
-VaR_{\alpha }\}} g((x-\mu )\Sigma^{-1}(x-\mu )^{t})  dx .
$$
Changing variables to $ y=(x-\mu)A^{-1} $ , $dy = \mid A\mid dx$ ,
where $\Sigma = A^{t} \; A $ is a Cholesky decomposition of $A $,
this becomes
$$
\alpha=  \int_{\{ \delta A \cdot y \leq -\delta \cdot \mu
-VaR_{\alpha } \}} g( | y |^{2} )  dy .
$$
Let $R $ be a rotation which sends $\delta A $ to $( |\delta A|,
0,...,0) $. Changing variables once more to $y=zR $, we obtain the
equation
$$
\alpha =  \int_{\{ |\delta A| z_{1} \leq - \delta \cdot \mu
-VaR_{\alpha } \}} g( | z |^{2} ) dz .
$$
If we write that $ | z |^{2} = z_{1}^2 + | z' |^{2} $ with $z' \in
\mathbb{R}^{n-1}$ then we have shown that :
$$
\alpha = \mbox{Prob } \{ \delta \cdot X < -VaR_{\alpha}\}=
\int_{\mathbf{R}^{n-1}} [ \int_{+\infty }^\frac{- \delta .\mu
-VaR_{\alpha }}{|\delta A| } g( z_{1}^2 + | z' |^{2}  ) dz_{1}]
dz' .
$$
Next, by using spherical  variables $z' = r\xi $ with $\xi\in
S_{n-2}$ , $dz' = r^{n-2} d\sigma(\xi ) dr $, we see that we have
to solve for $VaR_{\alpha } $ in the equation \be \alpha =
|S_{n-2}|  \int_{0}^{+\infty}  r^{n-2} \Big{[} \int_{-\infty
}^\frac{- \delta \mu^{t}   -VaR_{\alpha }}{|\delta A| } g( z_{1}^2
+ r^{2}  )  dz_{1}\Big{]} dr , \label{3} \ee $|S_{n-2}| $ being
the surface measure of the unit-sphere in $\mathbb{R }^{n - 1 } $:
$$
|S_{n-2}| = \frac{ 2 \pi^{\frac{n-1}{2}}}{\Gamma(\frac{n-1}{2})} .
$$
We now introduce the function
\begin{eqnarray}
G(s) &=& \frac{ 2 \pi^{\frac{n-1}{2}}}{\Gamma(\frac{n-1}{2})}
\int_{-\infty }^{-s } \Big{[} \int_{0}^{+\infty} r^{n-2} g(
z_{1}^2
+ r^{2} ) dr \Big{]} dz_{1} \nonumber \\
&=& \frac{ \pi^{\frac{n-1}{2}}}{\Gamma(\frac{n-1}{2})}
    \int_ s ^{-\infty }  \int_{z_{1}^2}^{+\infty}
    (u-z_{1}^2)^{\frac{n-3}{2}}   g(u) du dz_1
\label{3} ,
\end{eqnarray}
where for the second line we changed variables $u = r^2 + z_1 ^2
$. and replaced $z_1 $ by $-z_1 . $ We then have proved the
following result:

\begin{theorem} \label{VaR-elliptic}
Suppose that the portfolio's Profit \& Loss function over the time
window of interest is, to good approximation, given by $\Delta
\Pi=\delta_{1} X_{1}+\delta_{2} X_{2}+...+\delta_{n} X_{n}$, with
constant portfolio weights $\delta _j $. Suppose moreover that the
random vector $X = (X_1 , \cdots , X_n ) $ of underlying risk
factors follows a continuous elliptic distribution, with
probability density given by $f_X (x)= {|\Sigma|}^{-1} g((x-\mu
)\Sigma^{-1}(x-\mu )^{t}) $ where $\mu$ is the vector mean and
$\Sigma $ is the variance-covariance matrix, and where we suppose
that $g (s^2 ) $ is integrable over $\mathbb{R } $, continuous and
nowhere 0. Then the portfolio's {\em Delta-elliptic VaR} $VaR
_{\alpha } $ at confidence $1 - \alpha $ is given by
$$
VaR _{\alpha } =  - \delta \cdot \mu + q_{\alpha,n } ^g \cdot
\sqrt{ \delta \Sigma \delta ^t } ,
$$
where $q_{\alpha } = q_{\alpha,n } ^g $ is the, unique, solution
of the transcendental equation
$$
\alpha = G ( q_{\alpha,n } ) .
$$
\end{theorem}

\begin{remark} \rm{Note that $|\delta A | $ has a clear financial
interpretation, since \be |\delta A|=\sqrt{\delta \cdot \Sigma
\cdot \delta^t} , \label{rem1} \ee which is simply the portfolio's
volatility, or the square of its variance. }
\end{remark}

\begin{remark} \rm{In short-term Risk Management one can usually
assume that $\mu \simeq 0 $. In that case, we have that
$$VaR_{\alpha}=\sqrt{\delta \Sigma \delta^{t}} \cdot q_{\alpha,n } ^g $$
which is completely analogous to the result for linear portfolios
with normally distributed risk factors, except that, for example
for $\alpha = 0.05 $, the normal quantile at 5\%, which is
approximately 1.65, is now replaced by the $g $-dependent constant
$q^g _{0.05} $. The latter will have to be computed numerically,
for the different $g $'s one would like to use. }
\end{remark}

\begin{remark} \rm{One can in fact do the integral over $z_1 $ in
(\ref{3}): by Fubini,
\begin{equation}
G(s) = \int _s ^{\infty } K(s, u ) g(u) du , \label{3a}
\end{equation}
where the kernel $K $ is given by:
\begin{eqnarray*}
K(s, u ) &=& \frac{1 }{2 } |S_{n - 2 } | \int _{\sqrt{s} }
^{\sqrt{u } } (u - z_1 )^{\frac{n - 3 }{2 } } dz_1 \\
&=& \frac{1 }{4 } |S_{n - 2 } | \int _s ^u (u - y )^{\frac{n - 3
}{2 } } y^{-\frac{1 }{2 } } dy \\
&=& \frac{1 }{4 } |S_{n - 2 } | \int _0 ^{u - s } x^{\frac{n - 3
}{2 } } (u - x )^{-\frac{1 }{2 } } dx .
\end{eqnarray*}
At this stage we can the following integral from Gradshteyn and
Ryzhik \cite{GR}:
$$
\int _0 ^u \frac{x^{\mu - 1 } }{(1 + \beta x )^{\nu } } dx =
\frac{u^{\mu } }{\mu } \; _2 F_1 (\nu , \mu ; \mu + 1 ; - \beta u
) ,
$$
provided $\mbox{Re } \mu > 0 $ and $|\mbox{arg } (1 + \beta u ) |
< \pi $; cf. \cite{GR}, formula 3.194(1). It follows that
\begin{equation}
K(s, u ) = \frac{ \pi^{\frac{n-1}{2}}}{\Gamma(\frac{n+1}{2})} (u -
s )^{\frac{n - 1 }{2 } } \; _2 F_1 \left( \frac{1 }{2 } , \frac{n
- 1 }{2 } ; \frac{n + 1 }{2 } ; u(u - s ) \right) . \label{3b}
\end{equation}
However, we shall see in the example of the multi-variate $t
$-distribution which we will treat next, that it can be easier to
work directly with the double integral version (\ref{3}) instead
of with (\ref{3a}), (\ref{3b}). }
\end{remark}

\subsection{\bf{The case of t-student Distributions}}
\medskip
   We now consider in detail the case where our elliptic distribution is a multivariate
   Student-$t $. We will, unsurprisingly, call the corresponding $VaR$ the
   {\em Delta-Student VaR}, generalizing the familiar terminology
   of Delta Normal VaR.

In the case of multi-variate t-student  distributions, the
portfolio probability density function is given by:
$$
f_X (x) = \frac{\Gamma (\frac{\nu +
n}{2})}{\Gamma(\nu/2).\sqrt{|\Sigma|(\nu \pi)^n }}
{\Big{(}1+\frac{(x-\mu)^{t}\Sigma^{-1}(x-\mu) }{\nu}
\Big{)}}^{(\frac{-\nu-n}{2})} ,  $$ $ x \in \mathbb{R}^{n} $ and $
\nu > 2   $. Hence $g $ is given by
$$g(s)= C(\nu ,n) {(1+s/\nu )}^{-\frac{(n+\nu )}{2}} , \ \ s \geq 0 , $$
Where we have put
$$
C(\nu ,n)=\frac{\Gamma (\frac{\nu + n}{2})}{\Gamma(\nu/2)
\sqrt{(\nu \pi)^n }} .
$$
Using this $g $ in $(\ref{3})$, we find that

\be G (s) = \frac{\nu^\frac{n+\nu}{2} }{2 } |S_{n - 2 } | C(\nu ,
n )\int_s ^{\infty} I(z_1 ) dz_1 , \label{4} \ee

where we have put \be I(z_{1})=  \int_{z_{1}^2}^{+\infty}
(u-z_{1}^2)^{\frac{n-3}{2}} (\nu + u)^{-\frac{(n+\nu )}{2}}   du .
\label{5} \ee

The function $I(z_{1})$ can be evaluated with the help of another
one of the integrals in \cite{GR}:

\begin{lemma} \label{GR1} (Cf. \cite{GR} ,page 314.) If $|arg
(\frac{u}{\beta})|<\pi $, and $Re(\nu_{1})>Re(\mu)>0$ , then \be
\int_{w}^{+\infty}
    (x-w)^{\mu-1} (\beta + x)^{-\nu_{1}}
    dx={(w+\beta)}^{\mu-\nu_{1}} B(\nu_{1}-\mu,\mu) ,
    \label{rem2}
\ee with $B(\alpha , \beta ) $ the Euler Beta function:
$$
B(\alpha , \beta ) = \frac{\Gamma (\alpha ) \Gamma (\beta )
}{\Gamma (\alpha + \beta ) } .
$$
\end{lemma}

Using formula $(\ref{rem2})$ with $\nu_{1}=\frac{(n+\nu )}{2}$,
$\mu=\frac{n-1}{2}$,$\beta=\nu$, and $w = z_{1}^2$ , and
therefore, $\mu-\nu_{1}=-\frac{1+\nu}{2}$ and
$-\mu+\nu_{1}=\frac{1+\nu}{2}$, we find that

\be I(z_{1})={(z_{1}^2+\nu)}^{-\frac{1+\nu}{2}}
B(\frac{1+\nu}{2},\frac{n-1}{2})   \label{6}\ee

We have not finished yet, since we still have to integrate over
$z_1 $ in $(\ref{4})$. We therefore have to evaluate

\be J(s ,\nu)= \int_{-\infty }^{-s}
{(z_{1}^2+\nu)}^{-\frac{1+\nu}{2}} d{z_{1}} \label{7} \ee

Changing variable in this integral according to $u=z_{1}^2 $, we
find that
 \be  J(s ,\nu)= \frac{1}{2} \int_{s^2 } ^{\infty } u^{-\frac{1}{2}}
{(u+\nu)}^{-\frac{1+\nu}{2}} du \label{8} \ee

For the latter integral, we will use another formula from \cite{GR}:

\begin{lemma} (cf. \cite{GR}, formula 3.194(2)).
If $|arg (\frac{u}{\beta})|<\pi $, and $Re(\nu_{1})>Re(\mu)>0$ ,
then \be \int_{u}^{+\infty}
    x^{\mu-1} (1 + \beta x)^{-\nu_{1}}
    dx= \frac{u^{\mu-\nu_{1}}
     \beta^{-\nu_{1}}}{\nu_{1} -\mu }
    {_2F}_{1}(\nu_{1},\nu_{1}-\mu;\nu_{1}-\mu+1;-\frac{1}{\beta \cdot
    u}).
    \label{rem3}                  \ee
Here $_2 F _1 (\alpha ; \beta , \gamma ; w ) $ is the
hypergeometric function.
    \end{lemma}
In our case, $\nu_{1}=\frac{1+\nu}{2}$, $\mu=\frac{1}{2}$,
$\nu_{1}-\mu=\frac{\nu}{2}$,$\beta=\nu^{-1}$, and $u=s^2 $. If we
replace in $(\ref{8})$, we will obtain the following expression.
\be J(s ,\nu) = \frac{2}{\nu} s^{-\nu}
{_2F}_{1}\Big{(}\frac{1+\nu}{2},\frac{\nu}{2};1+\frac{\nu}{2};-\frac{\nu
}{s^2 } \Big{)}
 \label{9}  \ee
 Recalling (\ref{4}), we find, after a small computation, that in
 the Student-$\frak{t } $ case,
 \begin{eqnarray}
 G(s) = G_{\nu } ^{\frak{t } } (s) &=& \frac{1 }{\nu } \nu ^{\frac{n + \nu }{2 } } |S_{n - 2 } |
 C(\nu , n ) s^{- \nu }
 {_2F}_{1}\Big{(}\frac{1+\nu}{2},\frac{\nu}{2};1+\frac{\nu}{2};
 -\frac{\nu }{s^2 } \Big{)} \nonumber \\
 &=& \frac{1 }{\nu \sqrt{\pi } } \left( \frac{\nu }{s^2 } \right)
 ^{\nu / 2 } \frac{\Gamma \left( \frac{\nu + 1 }{2 } \right) }
 {\Gamma \left( \frac{\nu }{2 } \right) }
{_2F}_{1}\Big{(}\frac{1+\nu}{2},\frac{\nu}{2};1+\frac{\nu}{2};-\frac{\nu
}{s^2 } \Big{)} . \label{9a}
\end{eqnarray}
Hence we have proved the following result on Delta-Student VaR:

\begin{theorem}
Assuming that $\Delta \Pi  \simeq \delta_{1} X_{1}+\delta_{2}
X_{2}+...+\delta_{n} X_{n}$ with a multivariate Student-$\frak{t }
$ random vector $(X_{1},X_{2},..,X_{n})$ with vector mean $\mu$ ,
and variance-covariance matrix $\Sigma $, the linear Value-at-Risk
at confidence $1 - \alpha $ is given by the following formula
$$
VaR _{\alpha } = - \delta \cdot \mu +  q_{\alpha , \nu }^{\frak{t}
} \cdot \sqrt{ \delta \Sigma \delta ^t },
$$
where now $s = q_{\alpha }^t $ is the unique positive solution of
the transcendental equation
$$
G_{\nu } ^{\frak{t} } (q_{\alpha . \nu } ^{\frak{t} } ) = \alpha ,
$$
with $G_{\frak{t} } $ defined by (\ref{9a}).
\end{theorem}

\begin{remark} \rm{Note that $q_{\alpha , \nu } ^{\frak{t} } $ does not
depend of $n $. } \end{remark}

Hypergeometric $_2 F _1 $'s have been extensively studies, and
numerical software for their evaluation is available in Maple and
in Mathematica.
\\
\subsection{Some Numerical Result of Delta student VaR coefficient
$q_{\alpha,\nu}$ }

\medskip
   In the following table, we estimate only the positive solution of
   $G(s)=\alpha $ for some $\nu$ Given, with the help of Mathematica 4 Software.
$$
\begin{tabular}{|c|c|c|c|c|c|c|c|c|}

  \hline

$\nu$  & 2 & 3 & 4 & 5 & 6 & 7 & 8 & 9
  \\
\hline
  $q_{0.01,\nu}$   & 6.96456 & 4.54056  & 3.74695 & 3.36493 & 3.14267 & 2.99795  & 2.89646 & 2.8214  \\
 \hline
 $q_{0.025,\nu}$  &4.3026 & 3.18244  & 2.77644 & 2.57058 & 2.44691 & 2.36462  & 2.3060 & 2.26216    \\
 \hline
 $q_{0.05,\nu}$   & 2.91999 & 2.35336  & 2.13185 & 2.01505 & 1.94318 & 1.89458  & 1.85955 & 1.81246   \\
 \hline
\end{tabular}
$$
\\
\\
$$
\begin{tabular}{|c|c|c|c|c|c|c|c|c|}

  \hline

$\nu$  & 10 & 100 & 200 & 250 & 275 & 300 & 400 & 1000 \\
\hline
  $q_{0.01,\nu}$ & 2.76377 & 2.36422 & 2.34135 & 2.34514  & 2.33998 &  2.33884 & 2.33571 & 2.33008  \\
 \hline
  $q_{0.025,\nu}$ & 2.22814  &1.98397 & 1.97189  & 1.96949 & 1.96862 & 1.9679 & 1.96591 & 1.96234   \\
 \hline
$q_{0.05,\nu}$ & 1.66023  & 1.66023 & 1.65251 &  1.65097  & 1.65041 & 1.64995 & 1.64867 & 1.64638  \\
 \hline
\end{tabular}
$$
\begin{remark}  \rm{Note that, we obtain practically the same result as
in the case of normal distribution, when the degree of freedom of
our t-student is sufficiently high ($\nu $ near $300 $), as it of
course should. }
\end{remark}

\section{Linear VaR with mixtures of elliptic Distributions}
\medskip
\quad  Mixture distributions can be used to model situations where
the data can be viewed as arising from two or more distinct
classes of populations; see also \cite{MX}. For example, in the context
of Risk Management, if we divide trading days into two sets, quiet
days and hectic days, a mixture model will be based on the fact
that returns are moderate on quiet days, but can be unusually
large or small on hectic days. Practical applications of mixture
models to compute VaR can be found in Zangari (1996), who uses a
mixture normal to incorporate fat tails in VaR estimation. Here we
sketch how to generalize the preceding section to the situation
where the joint log-returns follow a mixture of elliptic
distributions, that is, a convex linear combination of elliptic
distributions.

\medskip
\begin{definition} \rm{We say that $ (X_{1},...,X_{n})$ has a joint distribution
that is the mixture  of $m $ elliptic distributions
$N(\mu_{j},\Sigma_{j},\phi_{j}) $\footnote{or
$N(\mu_{j},\Sigma_{j},g_{j})$ if we parameterize elliptical
distributions using $g $ instead of $\phi $}, with weights
$\{\beta_{j}\}$  (j=1,..,m ;  $\beta_{j} > 0$ ;  $\sum_{j=1}^m
\beta_{j} = 1$), if its cumulative distribution function can be
written as
$$
F_{X_{1},...,X_{n}}(x_{1},...,x_{n}) = \sum_{j=1}^m \beta_{j}
F_{j}(x_{1},...,x_{n})
$$
with $F_{j}(x_{1},...,x_{n}) $ the cdf of
$N(\mu_{j},\Sigma_{j},\phi_{j}) $. } \end{definition}

\begin{remark} \rm{
In practice, one would usually limit oneself to $m = 2 $, due to
estimation and identification problems; see \cite{MX}. }
\end{remark}

We will suppose that all our elliptic distributions
$N(\mu_{j},\Sigma_{j},\phi_{j}) $
admit a pdf :\\
\be f_{j}(x)= |\Sigma _j |^{-1/2 } g_{j}((x-\mu_{j} )
{\Sigma_{j}}^{-1}  (x-\mu_{j} )^t ). \ee The pdf of the mixture
will then simp;y be $\sum _{j = 1 } ^m \beta _j f_j (x) . $
\medskip

Let
$$\Sigma_{j} = A_{j} ^t \; A_{j} $$
be a Cholesky decomposition of $\Sigma _j . $ Since integration is
a linear operation, we now have to solve \be \alpha  = |S_{n-2}|
\sum_{j=1}^m \beta_{j}{|\Sigma_{j}|}^{-1/2} \int_{0}^{+\infty}
r^{n-2} \Big{[} \int_{-\infty }^\frac{- \delta \cdot \mu _j
-VaR_{\alpha }}{|\delta A_{j}| }      g_{j}( z_{1}^2 + r^{2}  )
dz_{1}\Big{]} dr \label{mix 1}\ee to obtain $VaR _{\alpha } $.
This leads to the following theorem:
\medskip
\begin{theorem}
Let $\Delta \Pi=\delta_{1} X_{1}+\ldots+\delta_{n} X_{n}$  with
 $(X_{1},\ldots,X_{n})$ is a mixture of elliptic distributions, with
 density
$$
f(x)=\sum_{j=1}^m  \beta _j {|\Sigma_{j}|}^{-1/2 }
g_{j}((x-\mu_{j} )\Sigma_{j}^{-1}(x-\mu_{j} )^{t})
$$
where $\mu_{j}$ is the vector mean, and  $\Sigma_{j} $ the
variance-covariance matrix of the $j $-th component of the
mixture.  We suppose that each $g_j $ is  integrable function over
$\mathbb{R}$, and that the $g_j $ never vanish jointly in a point
of $\mathbb{R }^m $. Then the value-at-Risk, or {\em Delta
mixture-elliptic VaR}, at confidence $1 - \alpha $ is given as the
solution of the transcendental equation
\begin{equation}
\alpha = \sum _{j = 1 } ^m \beta _j G_j \left( \frac{\delta
.\mu_{j}^{t} + VaR_{\alpha }}{(\delta \Sigma _j \delta )^{1/2 } }
\right) ,
\end{equation}
where $G_j $ is defined by (\ref{3}) with $g = g_j . $ Here
$\delta =(\delta_{1},\ldots,\delta_{n})$.
\end{theorem}
\begin{remark} \rm{
In the case of a mixture of m elliptic distributions the VaR will
not depend any more in a simple way on the total portfolio mean
and variance-covariance. This is unfortunate, but already the case
for a mixture of normal distributions. } \end{remark}
\begin{remark} \rm{
One might, in certain situations, try to model with a mixture of
elliptic distributions which all have the same variance-covariance
and the same mean, and obtain for example a mixture of different
tail behaviors by playing with the $g_j $'s. In that case the VaR
again simplifies to: $ VaR _{\alpha } = - \delta \cdot \mu +
q_{\alpha } \cdot \sqrt{ \delta \Sigma \delta ^t } $, with
$q_{\alpha } $ now the unique positive solution to
$$
\alpha  = \sum_{j=1}^m \beta_{j} G_j (q_{\alpha} ) .
$$
}
\end{remark}
\noindent The preceding can immediately be specialized to a
mixture of Student $\frak{t } $-distributions: the details will be
left to the reader.

\section{Expected Shortfall for elliptic distributions }
\medskip

Expected shortfall is a sub-additive risk statistic that describes
how large losses are on average when they exceed the VaR level.
Expected shortfall will therefore give an indication of the size
of extreme losses when the VaR threshold is breached. We will
evaluate the expected shortfall for a linear portfolio under the
hypothesis of elliptically distributed risk factors.
Mathematically, the expected shortfall associated with a given VaR
is defined as:
$$
\mbox{Expected Shortfall } = \mathbb{E } (-\Delta \Pi \vert
-\Delta \Pi > VaR ),
$$
see for example \cite{S}. Assuming again a multivariate elliptic
probability density $f(x) = {|\Sigma|}^{-1} g((x-\mu
)\Sigma^{-1}(x-\mu )^{t}) $, the Expected Shortfall at confidence
level $1 - \alpha $ is given by
\begin{eqnarray*}
- ES_{\alpha } &=& \mathbb{E } ( \Delta \Pi \mid   \Delta \Pi\leq
-VaR_{\alpha } ) \\
&=& \frac{1 }{\alpha } \mathbb{E } \left( \Delta \Pi \cdot 1_{\{
\Delta
\Pi \leq -VaR_{\alpha } \} } \right) \\
&=& \frac{1}{\alpha }\int_{\{ \delta x^t \leq -VaR_{\alpha }\}}
\delta x^t \ f(x) \ dx  \\
&=& \frac{ {|\Sigma|}^{-1/2}}{\alpha}  \int_{\{ \delta x^t \leq
-VaR_{\alpha }\}} \delta x^t \ g((x-\mu )\Sigma^{-1}(x-\mu )^{t})
dx .
\end{eqnarray*}
Let $\Sigma = A^t \; A $, as before.Doing the same linear changes
of variables as in section 2, we arrive at:
\begin{eqnarray*}
- ES_{\alpha} &=& \frac{1}{\alpha }  \int_{\{ |\delta A| z_{1}
\leq - \delta \cdot \mu - VaR_{\alpha } \}}  \ (|\delta A|z_{1} +
\delta \cdot \mu ) \   g( {\| z\|}^{2} )  dz \\
&=& \frac{1}{\alpha }  \int_{\{ |\delta A| z_{1} \leq - \delta
\cdot \mu - VaR_{\alpha } \}} \ |\delta A|z_{1} \ g( {\| z\|}^{2}
) \ dz \ + \ \delta \cdot \mu .
\end{eqnarray*}
The final integral on the right hand side can be treated as
before, by writing $ {\| z\|}^{2} = z_{1}^2 + {\| z^{'}\|}^{2} $
and introducing spherical coordinates $z^{'} = r\xi $, $\xi\in
S_{n-2} $, leading to:
$$
- ES_{\alpha } = \delta \cdot \mu  + \frac{|S_{n-2}|}{\alpha }
\int_{0}^{\infty} r^{n-2} \Big{[} \int_{   -\infty }^{\frac{-
\delta \mu^{t} -VaR_{\alpha }}{|\delta A| }}    {  |\delta A| \;
z_{1}} \   g( z_{1}^2 + r^{2} ) dz_{1}\Big{]} dr
$$
We now first change $z_1 $ into $-z_1 $, and then introduce $u =
z_1^2 + r^2 $, as before. If we recall that, by theorem
\ref{VaR-elliptic},
$$
q_{\alpha,n}^g = {\frac{ \delta \cdot \mu +VaR_{\alpha }}{|\delta
A| }}
$$
then, simply writing $q_{\alpha } $ for $q_{\alpha,n}^g $, we
arrive at:
\begin{eqnarray*}
ES_{\alpha } &=& - \delta \cdot \mu + |\delta A | \; \frac{|S_{n -
2 } | }{\alpha } \cdot \int _{q_{\alpha } } \int _{z_1 ^2 }
^{\infty } z_1 (u - z_1 ^2 )
^{\frac{n - 3 }{2 } } g(u) \ du \ dz_1 \\
&=& - \delta \cdot \mu + |\delta A | \; \frac{|S_{n - 2 } |
}{\alpha } \cdot \int _{q_{\alpha } ^2 } ^{\infty } \frac{1 }{n -
1 } \left( u - q_{\alpha } ^2 \right) ^{\frac{n - 1 }{2 } } \ g(u)
\ du ,
\end{eqnarray*}
since
$$
\int _{q_{\alpha } } ^{\sqrt{u } } z_1 \left( u - z_1 ^2 \right)
^{\frac{n - 3 }{2 } } = \frac{1 }{n - 1 } \left( u - q_{\alpha }
^2 \right) ^{\frac{n - 1 }{2 } } .
$$
After substituting the formula for $|S_{n - 2 } | $ and using the
functional equation for the $\Gamma $-function, $\Gamma (x + 1 ) =
x \Gamma (x) $, we arrive at the following result:

\begin{theorem} \label{EllipticES}
Suppose that the portfolio is linear in the risk-factors $X = (X_1
, \cdots , X_n )$: $\Delta \Pi=\delta \cdot X $ and that $X \sim
N(\mu ,\Sigma ,\phi ) $, with pdf $f(x)= {|\Sigma|}^{-1} g((x-\mu
)\Sigma^{-1}(x-\mu )^{t}) $. If we write $q_{\alpha } = $, then
the expected Shortfall
at level $\alpha $ is given by :\\
\be ES_{\alpha } = - \delta \cdot \mu + |\delta \Sigma \delta ^t
|^{1/2 } \cdot \frac{\pi^{\frac{n-1}{2}}}{\alpha \cdot
\Gamma(\frac{n + 1}{2})} \cdot \int _{(q_{\alpha , n } ^g )^2 }
^{\infty } \left( u - (q_{\alpha , n } ^g )^2 \right) ^{\frac{n -
1 }{2 } } \ g(u) \ du . \label{EESformula} \ee
\end{theorem}

\subsection{Application: Student Expected Shortfall}

\medskip
In the case of multi-variate t-student distributions we have that
$g(u)= C(\nu ,n) {(1+u/\nu )}^{-\frac{(n+\nu )}{2}} $, with $C(\nu
, n ) $ given in section 2. Let us momentarily write $q $ for
$q_{\alpha , \nu }^{\frak{t } } $. We can evaluate the integral in
(\ref{EESformula}) using lemma \ref{GR1}, as follows:
\begin{eqnarray*}
&\ & \int _{q^2 } ^{\infty } (u - q)^{\frac{n - 1 }{2 } } \left( 1
+ \frac{u
}{\nu } \right) ^{-\frac{n + \nu }{2 } } du \\
&=& \nu ^{\frac{n + \nu }{2 } } (q^2 + \nu )^{-(\frac{\nu - 1 }{2
} ) } B\left( \frac{\nu - 1 }{2 } , \frac{n + 1 }{ 2 } \right) .
\end{eqnarray*}
If we pose that :
$$es_{\alpha,\nu}= \frac{1 }{\alpha \cdot \sqrt{\pi }
} \frac{\Gamma \left( \frac{\nu - 1 }{2 } \right) }{\Gamma \left(
\frac{\nu }{2 } \right) } \nu ^{\nu / 2 } \left( (q_{\alpha , \nu
} ^{\frak{t } } )^2 + \nu \right) ^{-\left( \frac{\nu + 1 }{2 }
\right) }\label{es} $$
 After substitution in (\ref{EESformula}), we find, after some computations, the
following result:

\begin{theorem} The Expected Shortfall at confidence level $1 -
\alpha $ for a multi-variate Student-distributed linear portfolio
$\delta \cdot X $, with
$$
X \sim \frac{\Gamma (\frac{\nu + n}{2})}{\Gamma(\nu/2).
 \sqrt{|\Sigma|(\nu \pi)^n }}  {\Big{(}1+\frac{(x-\mu)^{t}
 \Sigma^{-1}(x-\mu) }{\nu} \Big{)}}^{-(\frac{\nu + n}{2})}  ,
$$
is given by:
\begin{eqnarray*}
ES_{\alpha , \nu } ^{\frak{t } } &=& - \delta \cdot \mu + |\delta
\Sigma \delta ^t |^{1/2 } \cdot \frac{1 }{\alpha \cdot \sqrt{\pi }
} \frac{\Gamma \left( \frac{\nu - 1 }{2 } \right) }{\Gamma \left(
\frac{\nu }{2 } \right) } \nu ^{\nu / 2 } \left( (q_{\alpha , \nu
} ^{\frak{t } } )^2 + \nu \right) ^{-\left( \frac{\nu + 1 }{2 }
\right) } \\
&=& - \delta \cdot \mu + |\delta \Sigma \delta ^t |^{1/2 } \cdot
\frac{1 }{\alpha \cdot \sqrt{\pi } } \frac{\Gamma \left( \frac{\nu
- 1 }{2 } \right) }{\Gamma \left( \frac{\nu }{2 } \right) } \nu
^{\nu / 2 } \left( \left( {\frac{ \delta \cdot \mu +VaR_{\alpha
}}{|\delta \Sigma \delta |^{1/2 } }} \right)^2 + \nu \right)
^{-\left( \frac{\nu + 1 }{2 } \right) } \\
&=& - \delta \cdot \mu +  es_{\alpha,\nu}\cdot |\delta \Sigma
\delta ^t |^{1/2 }
\end{eqnarray*}

\end{theorem}
\noindent The Expected Shortfall for a linear Student portfolio is
therefore given by a completely explicit formula, once the VaR is
known. Observe that, as for the VaR, the only dependence on the
portfolio dimension is through the portfolio mean $\delta \cdot
\mu $ and the portfolio variance $\delta \Sigma \delta ^t . $
\medskip

 With the help of
Matlab we obtain the the following table of values of $e_{\alpha
,\nu } $:
$$
\begin{tabular}{|c|c|c|c|c|c|c|c|c|}

  \hline

$\nu$  & 2 & 3 & 4 & 5 & 6 & 7 & 8 & 9
  \\
\hline
  $es_{0.01,\nu}$   & 5.5722 & 5.9309  & 5.7879 & 5.4555 & 5.0799 & 4.7160  & 4.3819 & 4.0818  \\
 \hline
 $es_{0.025,\nu}$  & 8.6113  &7.6777 &  6.8216 & 6.0676 & 5.4326 & 4.9032  & 4.4601 & 4.0862    \\
 \hline
 $es_{0.05,\nu}$   & 11.7123 & 9.0750  & 7.4966 & 6.3797 & 5.5457 & 4.9007  & 4.3880 &  3.9711  \\
 \hline
\end{tabular}
$$
\\
\\
$$
\begin{tabular}{|c|c|c|c|c|c|c|c|c|}

  \hline

$\nu$  & 10 & 100 & 200 & 250  \\ \hline
  $es_{0.01,\nu}$ & 3.8135 & 0.5157 & 0.2644 & 0.2086    \\
 \hline
  $es_{0.025,\nu}$ & 3.7675  & 0.4577 & 0.2313  & 0.1854    \\
 \hline
$es_{0.05,\nu}$ & 3.6257  & 0.4073 & 0.2050 &  0.1642    \\
 \hline
\end{tabular}
$$

\subsection{ES for Elliptic distribution Mixtures} The preceding
results can, as before, easily be generalized to mixtures of
elliptic distributions. Details will be left to the reader.

\section{Some Areas of Applications}
In this section, we survey some areas for applications of linear
portfolios that exist in the financial literature. We will discuss
4 examples:
\begin{itemize}
\item Delta-approximation of a derivatives portfolio. \item Linear
approximation of an equity portfolio \item Businesses as
portfolios of business units \item Incremental VaR
\end{itemize}

\subsection{Delta Approximation of a Portfolio}

\medskip

Suppose that we are holding a portfolio of derivatives depending
on $n $ underlying assets $X^{(1)},X^{(2)},...,X^{(n)}$with
elliptically distributed log-returns $r_j $ (over some fixed
time-window). The portfolio's present value $V $ will in general
be some complicated non-linear function of the $X_i $'s. To obtain
a first approximation of its VaR, we simply approximate the
present Value V of the position using a first order Taylor
expansion:
$$
V(X+\Delta X)\approx V(X) + \sum_{i=1}^n  \frac{\partial V}{
\partial X^{(i)}} \Delta X_i .
$$
>From this, we can then approximate the profit/loss function as
$$
\Delta V=V(X+\Delta X)- V(X)\approx \sum_{i=1}^n \delta_{i}
r^{(i)}=\delta \cdot r ,
$$
where we put $r=(r^{(1)},...,r^{(n)})$ and
$\delta=(\delta_{1},...,\delta_{n})$ with
$\delta_{i}=X^{(i)}\frac{\partial V}{ \partial X^{(i)}}$. The
entries of the $\delta$ vector are called the "delta equivalents "
for the position, and they can be interpreted as the sensitivities
of the position with respect to changes in each of the risk
factors. For more details see \cite{MX}, where a multi-variate normal
distribution for the $r_i $'s is assumed. The discussion there
generalizes straightforwardly to the elliptic case, where the
present paper's results can be used.

\subsection{Portfolios of Equities}
\medskip

A special case of the preceding is that of an equity portfolio,
build of stock $S_1 , \cdots , S_n $ with log-returns $r_1 (t) ,
\cdot , r_n (t) . $ In this case,
\begin{eqnarray*}
\Pi(t) - \Pi (0) &=& \sum_{i=1}^n w_{i} S_{i}(0) (
S_{i}(t)/S_{i}(0) -1 ) \\
&\approx & \sum _{i = 1 } ^n w_i S_i (0) r_i (t) ,
\end{eqnarray*}
where this approximation will be good if the $r_i (t) $ are small.

\subsection{Businesses as Linear Portfolios of Business Units}
An interesting way of looking upon an big enterprize, e.g. a
multi-national or a big financial institution, is by considering
it as a sum of its individual business units, cf. Dowd \cite{KD}. If
$X_{j}$, is the the variation of price or of profitability of
business unit j in one period , then the variation of price of the
agglomerate in the same period will be
$$
\Delta \Pi= X_{1}+\cdots+ X_{n} .
$$
The entire institution is therefore modelled by a linear
portfolio, with $\delta=(1, 1, \ldots , 1) $, to which the results
of this paper can be applied, if we model the vector of individual
price variations by a multi-variate elliptic distribution. VaR,
incremental VaR (see below) and Expected Shortfall will be
relevant here. For more details see Dowd \cite{KD}, chapter XI .

\subsection{Incremental VaR}
Incremental VaR is defined in \cite{MX}  as the statistic that provides
information regarding the  sensitivity of VaR to changes in the
portfolio holdings. It therefore gives an estimation of the change
in VaR resulting from a risk management decision. Results from
\cite{S} for incremental VaR with normally distributed
risk-factors generalize straightforwardly to elliptically
distributed ones: if we denote by $IVaR_{i}$ the incremental VaR
for each position in the portfolio, with $\theta_{i}$ the
percentage change in size of each position,then the change in VaR
will be given by
$$
\Delta VaR=\sum \theta_{i}IVaR_{i} \label{ivar1}
$$
By using the definition of $IVaR_{i}$ as in \cite{MX}  (2001), we have
that \be IVaR_{i}=\omega_{i}\frac{\partial VaR }{\partial
\omega_{i}} \ee with $\omega_{i}$ is the amount of money invested
in instrument $i $. In the case of an equity portfolio in the
elliptically distributed assets, we have seen that, assuming
$\mu=0$,
$$
VaR_{\alpha} = - q_{\alpha,n}^g \sqrt{\delta \Sigma \delta^t } ,
$$
We can then calculate $IVaR_{i}$ for the i-th constituent of
portfolio as
$$
IVaR_{i}=\omega_{i}\frac{\partial VaR }{\partial
\omega_{i}}=\omega_{i}\gamma_{i}
$$
with
$$
\gamma= - q_{\alpha,n}^g \frac{\Sigma \omega}{\sqrt{\delta \Sigma
\delta^t } } .
$$
The vector $\gamma $ can be interpreted as a gradient of
sensitivities of VaR with respect to the risk factors. This is the
same as in \cite{MX}, except of course that the quantile has changed
from the normal one to the one associated to $g . $

\section{conclusion}
\medskip
In this paper we have shown how to reduce the estimation of
Value-at-Risk for linear elliptic portfolios to the evaluation of
one dimensional integrals which, for the special case of a Student
distribution, can be explicitly evaluated in terms of a
hypergeometric function. We indicated how to extend these to the
case of mixtures of elliptic distributions. We have also given a
similar, but simpler, integral formula for the expected shortfall
of such portfolios which, again, can be completely evaluated in
the Student case. We finally surveyed some potential application
areas.


\newpage


\baselineskip=10pt
\begin{thebibliography}{99}
\bibitem{AS} Albenese, C., and L.Seco :  "Harmonic Analysis in Value-at-Risk Calculations" , { \it  Revista Matemaatica Ibero Americana 17 (2),2001}.
\bibitem{DP}  Duffie, D., and  J.Pan :  An overview of Value at Risk , {\it J. Derivatives 4(3), 7-49.}
\bibitem{EMS}  P.Embrechts , A.McNeil and D.Strauman (1999)
Correlation and dependance in Risk Management      Properties and
Pitfalls. In M.A.H. Dempster, editor, Risk Management : Value at Risk and Beyond, pages 176-223. Cambridge University Press, 1999.
\bibitem{GR} {I.S Gradshteyn, I.M. Ryzhik : Table of integrals, series, and products  ( 2000)   Editor Alan Jeffrey }
\bibitem{GHS}Glasserman, P. , Heidelberger, P. and Shahabuddin :
 Portfolio Value-at-Risk with Heavy Yailed Risk Factors. { \it Mathematical Finance, 12 (3): 239-269, July 2002}.
\bibitem{MX} Jorge Mina and Jerry Yi Xiao (2001), Return to
Riskmetrics: The Evolution of the standard (www.riskmetrics.com).
2001 RiskMetrics Group Inc.
\bibitem{KD}  Kevin Dowd,(1998) Beyond Value-at-Risk: the new science
of risk Management (Wiley series in Frontiers in Finance).
\bibitem{LW} Lopez, J.A. and Walter,C.A.(2000). Evaluating
covariance matrix Forecastsin a Value-at-Risk Framework, FRBSF
working paper 2000-21, Federal Reserve Bank of San Francisco
\bibitem{KB} Kotz, Balahskrinan, (2001)    Continuous Multivariate Distributions
\bibitem{S} J.Sadefo-Kamdem ,(2003):VaR Estimation for  Portfolio of Securities with joint elliptic distribution  log-Returns.    Working paper
\bibitem{YT} Yasuhiro and Toshinao Yoshiba, (2002): On the validity of Value-at-Risk: Comparative analysis with Expected shortfall
\bibitem{Z} Zangari, P. (1996). An Improved Methodology for
Measuring VaR, RiskMetrics Monitor, 2nd quarter pp.7-25.
http://www.riskmetrics.com/research/journals
\end{thebibliography}
\end{document}